\documentclass{article}
\usepackage{nccmath}
\usepackage[mathscr]{euscript}
\usepackage{tikz}
\usepackage{dsfont}
\usepackage{color}
\usepackage{subfig}
\usepackage{xcolor}

\numberwithin{equation}{section}

\definecolor{light}{gray}{.75}
\newcommand{\R}{\mathds{R}}

\newtheorem{theorem}{Theorem}[section]

\newtheorem{definition}{Definition}[section]
\newtheorem{example}{Example}[section]

\def\clap#1{\hbox to 0pt{\hss#1\hss}}

\makeindex             % used for the subject index
\begin{document}
\title{Schur functions for approximation problems}
\author{Nadezda Sukhorukova and Julien Ugon}
%\institute
%{A. Fontes \at Faculty of Science, Engineering and Technology, Swinburne University of Technology, PO Box~218, Hawthorn, Victoria, Australia, \email{afontes@student.swin.edu.au}
%\and 
%{Nadezda Sukhorukova \at Faculty of Science, Engineering and Technology, Swinburne University of Technology, PO Box~218, Hawthorn, Victoria, Australia and  Centre for Informatics and Applied Optimization, Federation University Australia \email{nsukhorukova@swin.edu.au}
%and
%Julien Ugon \at Faculty of Science, Engineering and Built Environment, Deakin University, 221 Burwood Highway
%Burwood Victoria 3125, Australia and  Centre for Informatics and Applied Optimization, Federation University Australia \email{julien.ugon@deakin.edu.au}
\maketitle

\abstract{
%In this paper we propose a new efficient algorithm for signal clustering. This approach is based on a $k$-means type clustering and %least squares approximation. There are two main computational problems here. Firstly, we need to find a least squares %approximation (prototype) of a group of curves (curve cluster).  Secondly, we need to develop a fast and accurate procedure for %cluster prototype updates when new segments of signal for the same time period become available. The results of numerical %experiments are provided.
In this paper we propose a new approach to least squares approximation problems. This approach is based on partitioning and Schur function. The nature of this approach is combinatorial, while most existing approaches are based on algebra and algebraic geometry. This problem has several practical applications. One of them is curve clustering. We use this application to illustrate the results.}

\section{Introduction}
\label{sec:intro}

In this paper we formulate a specific least  least squares approximation problem and provide  a signal processing application where this problem is used.   The main technical difficulty for this problem is to solve linear systems with the same system matrix and different right-hand sides. One simple approach that can be proposed here is to invert the system matrix and multiply the updated right-hand side by this inverse at each iteration. In general, it is not very efficient to solve linear systems through computing matrix inverses, but in this particular application it is very beneficial. One technical difficulty here is to know in advance whether the system matrix is invertible or not. Similar problems appear in Chebyshev (uniform) approximation problems as well.

 In this paper we suggest a new approach for dealing with this kind of systems. This approach in based on Schur functions, a well-established techniques that is used to describe partitioning~\cite{macdonald1995symmetric}. The very nature of these functions is  combinatorial. Based on our previous experience~\cite{MATRIX2016}, the characterisation of the necessary and sufficient optimality conditions for multivariate Chebyshev approximation is also combinatorial and therefore Schur function is a very natural tool to work with these problems.

 This paper is organised as follows. In section~\ref{sec:application} we introduce a signal processing application that relies on approximation and optimisation. In section~\ref{sec:math_formulation} we provide a  mathematical formulation to the signal processing problem and discuss how it can be simplified. In section~\ref{sec:Schur} we introduce an innovative approach for solving the problem. This approach is based on Schur functions. Finally, in section~\ref{sec:discussion} we provide future research directions.

\section{Signal clustering}\label{sec:application}

In signal processing, there is a need for constructing signal prototypes. Signal prototypes are summary curves that may replace the whole group of signal segments, where the signals  are believed to be similar to each other. Signal prototypes may be used  for characterising the structure of the signal segments and also for reducing the amount of information to be stored. 

Any signal group prototype should be an accurate approximation for each member of the group. On the top of this, it is desirable that the process of recomputing group prototypes, when new group members are available, is not computationally expensive. 

 In this paper we suggest a $k$-means and least square approximation based model. Similar models are proposed in~\cite{spath80}.  This is a convex optimisation problem. There are several advantages of this model. First of all, it  provides an accurate approximation to the group of signals. Second, this problem can be obtained as a solution to a linear system  and can be solved efficiently. Finally, the proposed  approach  allows one to compute prototype updates without recomputing  from  scratch.

\section{Mathematical formulation}\label{sec:math_formulation}
\subsection{Prototype construction}
Assume that there is a group of $l$ signals $S_1(t),\dots,S_l(t)$, whose values are measured at discrete time moments $$t_1,\dots, t_N,~t_i\in[a,b],~i=1,\dots,N.$$
We suggest to construct the prototype as a polynomial $P_n({\bf X},t)=\sum_{i=0}^{n}x_it^i$ of degree $n$, whose least squares  deviation from each member of the group on $[a,b]$ is minimal. That is, one has to solve the following optimisation problem:
\begin{equation}\label{eq:main_pol_convex}
{\rm minimise}~ F({\bf X})= \sum_{i=1}^{N}\sum_{j=1}^{l}(S_j(t_i)-P_n({\bf X},t_i))^2,
\end{equation} 
where ${\bf X}=(x_0,\dots,x_n)\in\R^{n+1}$, $x_k,~k=0,\dots,n$ are the polynomial parameter and also the decision variables. Each signal is a column vector 
$${\bf S}^j=(S_j(t_1),\dots,S_j(t_N))^T,~j=1,\dots, l.$$

 Problem~(\ref{eq:main_pol_convex})  can be formulated in the following matrix form:
\begin{equation}\label{eq:main_pol_convex_matrix}
{\rm minimise}~ F({\bf X})= \|{\bf Y}-{\bf BX}\|,
\end{equation} 
where 
\begin{description}
\item 
${\bf X}=(x_0,\dots,x_n)\in\R^{n+1}$,  are the decision variables (same as in~(\ref{eq:main_pol_convex}));
\item  vector $${\bf Y}=\left(
\begin{matrix}{\bf S}^1\\
{\bf S}^2\\
\vdots\\
{\bf S}^l\end{matrix}
\right)\in\R^{(n+1)l}$$
\item    matrix ${\bf B}$ contains repeated matrix blocks, namely,
$${\bf B}=\left(
\begin{matrix}
{\bf B}_0\\
{\bf B}_0\\
{\bf B}_0\\
\vdots\\
{\bf B}_0\\
\end{matrix}
\right),$$
    where 
    $$
    {\bf B}_0=
    \left(
    \begin{array}{ccccc}
    1&t_1&t_1^2&\dots&t_1^n\\
     1&t_2&t_2^2&\dots&t_2^n\\
      \vdots&\vdots&\ddots&\dots&\vdots\\
       1&t_N&t_N^2&\dots&t_N^n\\
    \end{array}
    \right).
    $$
    \end{description}
This least squares problem can be solved using a system of normal equations: 
\begin{equation}\label{eq:normal_equations}
{\bf B}^T{\bf B}{\bf X}={\bf B}^T{\bf Y}.
\end{equation}
Taking into account the structure of  the system matrix in~(\ref{eq:normal_equations}), the problem can be significantly simplified:
\begin{equation}\label{eq:normal_equations1}
l{\bf B}_0^T{\bf B}_0{\bf X}={\bf B}_0^T\sum_{k=1}^{l}{\bf S}^k.
\end{equation}
Therefore, instead of solving~(\ref{eq:normal_equations}), one can solve
\begin{equation}\label{eq:normal_equations2}
{\bf B}_0^T{\bf B}_0{\bf X}={\bf B}_0^T{\sum_{k=1}^{l}{\bf S}^k\over l}={\bf B}_0^T{\bf S},
\end{equation}
where {\bf S} is the average of all $l$ signals of the group (centroid).

\subsection{Prototype update}
Suppose that a signal group prototype  has been constructed. 
%and $$({\bf X}^*)=(x^*_0,\dots,x^*_{n})$$ is the corresponding optimal solution to~(\ref{eq:normal_equations2}). 
Assume now that we need to update our group of signals: some new signals have to be included, while some others are to be excluded.
 To update the prototype, one needs to update the centroid and  solve~(\ref{eq:normal_equations2}) with the updated right-hand side, while the system matrix ${\bf B}_0^T{\bf B}$ remains the same. 
 
 If only few signals are moving in  and out of the group, then the updated centroid can be calculated without recomputing from scratch. Assume that $l_a$ signals are moving in the group (signals $S_a^1(t),\dots S_a^{l_a}$), while $l_r$ are moving out (signals $S_r^1(t),\dots S_r^{l_r}$), then the centroid can be recalculated as follows:
 $$S_{new}(t)={lS_{old}(t)+\sum_{k=1}^{l_a} S_{a}^k(t)+\sum_{k=1}^{l_r}S_r^k(t)\over l-l_r+l_a}.$$
 
Since the same system has to be solved repeatedly with different right-hand sides, one approach is to invert matrix~${\bf B}_0^T{\bf B}_0$, which is an $(n+1)\times(n+1)$ matrix. In most cases, $n$ is much smaller than $N$ or $l$ and therefore this approach is quite attractive, if we can guarantee that  matrix~${\bf B}_0^T{\bf B}_0$  is invertible. In the next section we discuss the verification of this property.

\section{Schur functions and matrix inverse}\label{sec:Schur}
\subsection{Vandermonde and generalised Vandermonde matrices}
Consider matrix~${\bf B}_0^T{\bf B}_0$. 	In general, matrix~${\bf B}_0$ can be defined as follows:
$$
    {\bf B}_0=
    \left(
    \begin{array}{ccccc}
    g_1(t_1)&g_2(t_1)&g_3(t_1)&\dots&g_{n+1}(t_1)\\
       g_1(t_2)&g_2(t_2)&g_3(t_2)&\dots&g_{n+1}(t_2)\\
      \vdots&\vdots&\ddots&\dots&\vdots\\
      g_1(t_N)&g_2(t_N)&g_3(t_N)&\dots&g_{n+1}(t_N)\\
    \end{array}
    \right),
    $$
    where $g_i,~i=1,\dots,n+1$ are basis functions.  In section~\ref{sec:math_formulation} we were discussing polynomial approximation and therefore, the components of matrix~${\bf B}_0^T{\bf B}_0$ are monomials that are evaluated at different time-moments. Recall that $n+1<<N$. Matrix~${\bf B}_0^T{\bf B}_0$ is invertible if and only if  matrix~${\bf B}_0$ has exactly $n+1$ linearly independent rows. This is always the case when functions $g_i,~i=1,\dots,n+1$ form a Chebyshev system (for example, monomials $g_i=t^{i-1}$, some systems of trigonometric functions). This is not always the case when, for example, some of the monomials are ``missing'' from the system. This situation is illustrated in the following example.
\begin{example}
Consider the system of two monomials on the segment~$[-1,1]$:
$$g_1(t)=1,~g_2(t)=t^2, ~t_1\neq t_2,$$
the monomial $t$ is ``missing''. Take time-moments $t_1, t_2\in[-1,1]$. The determinant
$$\left|
    \begin{array}{cc}
    1&t_1^2\\
     1&t_2^2\\
         \end{array}
    \right|=0\Leftrightarrow t_1=-t_2.
    $$
    Therefore, these functions do not form a Chebyshev system, since the corresponding determinant is zero when, for example, $t_1=-t_2=1$ and there is only one linear independent row. 
\end{example}

Recall that in the case of classical polynomial approximation (all monomials are included into the set of basis functions), the corresponding determinant is non-zero as it is the determinant of a Vandermonde matrix. We now need to introduce so called generalised Vandermonde matrices.
\begin{definition}
Generalised Vandermonde matrices have the following structure:
$$G= \left(
\begin{matrix}
t_1^{m_1}& t_2^{m_1}&\dots &t_{n+1}^{m_1}\\
t_1^{m_2}& t_2^{m_2}&\dots &t_{n+1}^{m_2}\\
\vdots& \vdots&\ddots &\vdots\\
t_1^{m_n}& t_2^{m_n}&\dots &t_{n+1}^{m_n}\\
\end{matrix}\right).
$$
\end{definition}

%Assume that $0\leq m_n<m_{n-1}<\dots <m_1$.

Denote 
\begin{equation}\label{eq:degrees}
m_1=\lambda_1+n-1, ~m_2=\lambda_2+n-2,\dots, m_n=\lambda_n+n-n=\lambda_n.
\end{equation}
Define the following function
\begin{equation}\label{eq:schur_fun}
s_{\lambda}(t_1,\dots,t_{n+1})={{\rm det}(G)\over{\rm det} (V)},
\end{equation}
where $G$ is the matrix with one or more missing monomials  and $V$ is the Vandermonde matrix. Vandermonde matrices correspond to $$\lambda_1=\lambda_2=\dots=\lambda_n=0.$$

$s_{\lambda}(t_1,\dots,t_{n+1})$ is called Schur function, named after Issai Schur. Schur polynomials are certain symmetric polynomials of $n$ variables. This polynomials are used in representation and partitioning. A good introduction to Schur polynomials can be found in~\cite{macdonald1995symmetric}.  Therefore,
$${\rm det}(G)=s_{\lambda}(t_1,\dots,t_n){\rm det} (V)$$
and hence one needs to study the behaivour of Schur functions. 
Therefore, the following theorem holds.
\begin{theorem}\label{thm:schur_nonsingular}
Matrix  ${\bf B}^T_0{\bf B}_0$ is non-singular if and only if the corresponding Schur function~(\ref{eq:schur_fun}) is non-zero.
\end{theorem}

In particular, if $t_i>0,~i=1,\dots,n+1$, then the system is Chebyshev. Note that this statement can  be proven using a logarithmic transformation~\cite{Karlin66}. We believe, however, that our approach is also applicable to more general settings.  

There are many studies on Schur polynomials and many efficient ways for computing them.  This approach can be used, for example, if one needs to know if ${\bf B}^T_0{\bf B}_0$ is invertible. If the matrix is invertible, one can develop a very fast and efficient algorithm for curve cluster prototype updates.  If the matrix is singular, one can use the singular-value decomposition for constructing the prototype updates. This decomposition can be computed once, since ${\bf B}_0^T{\bf B}_0$ remains unchanged when the cluster membership is updated. 

\section{Discussions and future research directions}\label{sec:discussion}

There are many studies on how to compute Schur functions. We are particularly interested in the extension of this approach to Chebyshev (uniform) approximation and  multivariate approximation. This is a very promising approach for dealing with this type of problems, since, as our previous studies suggested~\cite{MATRIX2016} the corresponding optimality conditions are very combinatorial in their nature and therefore, Schur functions are a very natural tool for study this kind of systems.

We are also planning to conduct a thorough numerical study of the signal processing application we are discussing in this paper. 

\section{Acknowledgement}
This paper was inspired by the discussions during a recent
MATRIX program ``Approximation  Optimisation and Algebraic Geometry\rq{}\rq{}  that took place in February 2018.  We are thankful to the MATRIX organisers,  support team and participants for a terrific research atmosphere and productive discussions.

%    \bibliographystyle{amsplain}
 %   \bibliography{mybib1}
\end{document}